\documentclass[11pt]{amsart}
\usepackage{amssymb}
\usepackage{amsmath}
\usepackage{bbm}
\usepackage{graphicx}
\usepackage{hyperref}

\numberwithin{equation}{section}

\newcommand{\be}{\begin{equation}}
\newcommand{\ee}{\end{equation}}
\newcommand{\ds}{\displaystyle}
\newcommand{\ts}{\textstyle}

\newcommand{\bmubar}{\boldsymbol{\bar\mu}}
\newcommand{\wM}{{}^\nabla\! M}
\newcommand{\wB}{{}^\nabla\! B}

\DeclareMathOperator{\diag}{diag}

\newtheorem{thm}{Theorem}
\newtheorem{cor}[thm]{Corollary}
\newtheorem{lem}[thm]{Lemma}
\newtheorem{prp}[thm]{Proposition}
\newtheorem{defn}{Definition}

\makeatletter
\def\Ddots{\mathinner{\mkern1mu\raise\p@
\vbox{\kern7\p@\hbox{.}}\mkern2mu
\raise4\p@\hbox{.}\mkern2mu\raise7\p@\hbox{.}\mkern1mu}}
\makeatother

\title{Stationary Distribution and Eigenvalues \\ for a de Bruijn Process}
\author{Arvind Ayyer}
\address{Department of Mathematics \\
University of California \\
Davis CA 95616}
\email{ayyer@math.ucdavis.edu}
\author{Volker Strehl} 
\address{
Department of Computer Science \\ Universit\"at Erlangen-N\"urnberg \\
D-91058 Erlangen \\ Germany }
\email{strehl@informatik.uni-erlangen.de}

\dedicatory{Dedicated to Herb Wilf on the occasion of his 80th birthday}

\date{\today}

\begin{document}

\begin{abstract}
We define a de Bruijn process with parameters $n$ and $L$ as a certain 
con\-ti\-nu\-ous-time Markov chain on the de Bruijn graph with words of length $L$ 
over an $n$-letter alphabet as vertices. We determine explicitly its steady state 
distribution and its characteristic polynomial, which turns out to decompose 
into linear factors. In addition, we examine the stationary state of two 
specializations in detail. In the first one, the de Bruijn-Bernoulli process, this is 
a product measure. In the second one, the Skin-deep de Bruin process, the 
distribution has constant density but nontrivial correlation functions. 
The two point correlation function is determined using generating function techniques.
\end{abstract}

\maketitle

\section{Introduction}

A de Bruijn sequence (or cycle) over an alphabet of $n$ letters and of order $L$ 
is a cyclic word of length $n^{L}$ such that every possible word of length $L$ over 
the alphabet appears once and exactly once. The existence of such sequences and 
their counting was first given by Camille Flye Sainte-Marie in 1894 for the case $n=2$,
see \cite{flye} and the acknowledgement by de Bruijn\cite{bruijn75}, although the 
earliest known example comes from the Sanskrit prosodist Pingala's {\em Chandah 
Shaastra} (some time between the second century BCE and the fourth century CE 
\cite{vannooten,knuth4}). 
This example is for $n=2$ and $L=3$  essentially contains the word 0111010001 
as a mnemonic for a rule in Sanskrit grammar. Omitting the last two letters (since 
they are repeating the first two) gives a de Bruijn cycle.
Methods for constructing de Bruijn cycles are discussed by Knuth \cite{knuth2}.

The number of de Bruijn cycles for  alphabet size $n=2$ was (re-)proven to be 
$2^{2^{L-1}-L}$ de Bruijn  \cite{bruijn}, hence the name.  
The generalization to arbitrary alphabet size $n$ was first proven to be 
 $\ds n!^{n^{L-1}} \cdot n^{-L}$ by de Bruijn and van Aardenne-Ehrenfest. 
 This result can be seen as an application of the famous 
 BEST-theorem \cite{vanbruijn,smithtutte,tutte}, which relates the counting
 of Eulerian tours in digraphs to the evaluation of a Kirchhoff (spanning-tree counting)
 determinant. The relevant determinant evaluation for the case
 of de Bruijn graphs (see below) is due to Dawson and Good \cite{dawsongood}, 
 see also \cite{knuth67} .

The (directed) de Bruijn graph $G^{n,L}$ is defined over an alphabet $\Sigma$ of cardinality $n$.
Its vertices are the words of $u=u_1u_2\ldots u^L \in \sigma^L$, and there is an directed edge
or arc between any two nodes $u=u_1u_2\ldots u_L$ and $v=v_1v_2\ldots v_L$ if an only if
$t(u)=u_2 \ldots u_n = v_1 \ldots v_{n-1}=h(v)$, where $h(v)$ ($t(u)$ resp.) stands for
the \emph{head} of $v$ (\emph{tail} of $u$, resp.). 
This arc is naturally labeled by the word $w=u.v_L=u_1.v$,
so that $h(w)=u$ and $t(v)=v$. It is intuitively clear that Eulerian tours in the de Bruijn graph
$G^{n,L}$ correspond to de Bruijn cycles for words over $\Sigma$ of length $L+1$ 
de Bruijn graphs and cycles have applications in several fields, e.g. in
networking \cite{kaashkarg} and bioinformatics \cite{pevzneretal}.	
 For an introduction to de Bruijn graphs, see e.g. \cite{ralston}.
 
 In this article we will study a natural continuous-time Markov chain on $G^{n,L}$
which exhibits a very rich algebraic structure. The transition probabilities are
not uniform since they depend on the structure of the vertices as words,
and they are symbolic in the sense that variables are attached to the edges
as weights.
We have not found this in the literature, although there are studies of the 
uniform random walk on the de Bruijn graph \cite{flajoletetal}.
The hitting times \cite{chenh} and covering times \cite{mori} of this random walk 
have been studied, as has the structure of the covariance matrix for 
the alphabet of size $n=2$ \cite{alhakim-molchanov} and in general \cite{alhakim1}.
The spectrum for the undirected de Bruijn graph has been found by Strok \cite{strok}.
We have also found a similar Markov chain whose spectrum is completely 
determined in the context of cryptography \cite{geiselgoll}.

After describing our model on $G^{n,L}$ for a de Bruijn process in detail in the next section,
we will determine its stationary distribution in Section 3 and its spectrum in Section 4.
In the last section we discuss two special cases, the de Bruijn-Bernoulli process 
and the Skin-deep de Bruijn process.

\section{The Model}

We take the de Bruijn graph $G^{n,L}$ as defined above. As alphabet we may
take $\Sigma=\Sigma_n= \{1,2, \ldots,n\}$. Matrices will then be indexed by words over 
$\Sigma_n$ taken in lexicographical order. Since the alphabet size $n$ will be 
fixed throughout the article, we will occasionally drop $n$ as super- or subscript if there is no 
danger of ambiguity.

From each vertex $u=u_1u_2 \ldots u_L \in \Sigma^L$ there are $n$ directed edges 
in $G^{n,L}$ joining $u$ with the vertices $u_2u_3 \ldots u_n.a = t(u).a$ for $a \in \Sigma$.

We now give weights to the edges of the graph $G^{n,L}$.
Let $X =  \{x_{a,k} \,;\, a \in \Sigma, k \geq 1\}$ be the set of weights, 
to be thought of as formal variables.
An {\em $a$-block} is a word $u \in \Sigma^+ $ 
which is the repetition of the single letter $a$ so that
$u=a^{k}$ for some $a \in \Sigma$ and $k \geq 1$.
Obviously, 
every word $u$ has a unique decomposition into blocks of maximal length, 
\be \label{blockfac}
u = b^{(1)} b^{(2)} \cdots b^{(m)},
\ee
where each factor $b^{(i)}$ is a block so that any two neighboring factors  
are blocks of \emph{distinct} letters. This is the canonical block factorization of $u$ 
with a minimum number of block-factors.

We now define the function $\beta: \Sigma^{+} \to X$ as follows:
\begin{itemize}
\item[--] for a block $a^k$ we set $\beta(a^k)=x_{a,k}$; 
\item[--] for $u \in \Sigma^+$ with canonical block factorization (\ref{blockfac})
	we set $\beta(u) = \beta(u^{(m)})$, \\ i.e.,
	the $\beta$-value of the last block of $u$.
\end{itemize}
An edge from vertex $u\in \Sigma^L$ to vertex $v\in \Sigma^L$, 
so that $h(v)=t(u)$ with $v=t(u).a$, say,
will then be given the weight $\beta(v)$. This means that
\be
\beta(v) =
\begin{cases}   x_{a,L} & \text{if}~\beta(u)=x_{a,L}, \cr
			x_{a,k+1}&\text{if}~\beta(u)=x_{a,k}~\text{with}~k<L, \cr
                       x_{a,1}      &\text{if}~\beta(u)=x_{b,k}~\text{for some}~b \neq a.
\end{cases}
\ee

Our \emph{de Bruijn process} will be a  continuous time Markov chain
derived from the Markov chain represented by the directed de Bruijn graph $G^{n,L}$
with edge weights as defined above. The transition rates are $\beta(v)$ for transitions 
represented by edges ending in $v$. 
We note that these rates can be taken just as variables and not necessarily probabilities.
Similarly, expectation values of random variables in this process will be
functions in these variables.

\begin{figure}[h]
\begin{center}
 \includegraphics[width=5cm]{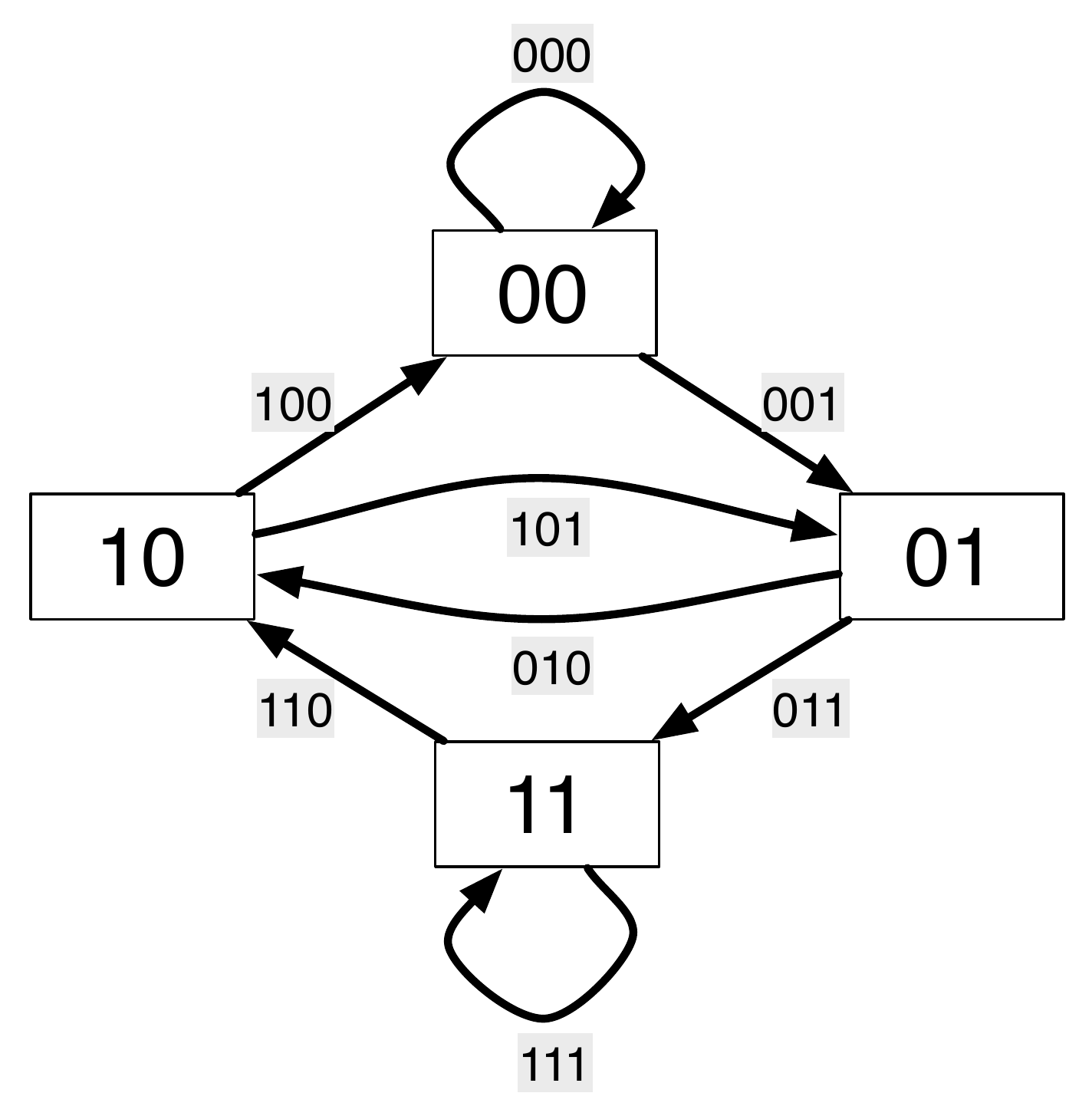} 
\caption{An example of a de Bruijn graph in two letters and words of length 2.}
\label{fig:example}
\end{center}
\end{figure}

The simplest nontrivial example occurs when $n=L=2$. There are four
 configurations and the relevant edges are given in the Figure~\ref{fig:example}. 
 
Before stating our notation for the transition matrix of a continuous-time Markov chain, 
our \emph{de Bruijn process}, we need a general notion.
\begin{defn} \label{def:tilde}
For any $k \times k$
 matrix $M$, let  ${}^\nabla\!M$ denote the matrix where the sum of 
each column is subtracted from the corresponding diagonal element,
	\be
	{}^\nabla \! M = M - \diag(1_{k} \cdot M),
	\ee
where $1_{k}$ denotes the all-one row vector of length $n$ and 
$\diag(m_{1}, \dots,m_{k})$ is a diagonal matrix with entries 
$m_{1}, \dots, m_{k}$ on the diagonal.
\end{defn}

In graph theoretic terms ${}^\nabla\!M$ is the (negative of) 
the \emph{Kirchhoff} matrix  or {\em Laplacian matrix} of G,
if $M$ is the weighted adjacency matrix of a directed graph $G$. 
In case $M$ is a matrix representing transitions of a Markov chain, 
the column (or right) eigenvector of ${}^\nabla\!M$
for eigenvalue zero properly normalized
gives the stationary probability distribution of
the continuous-time Markov chain.

We note that the graphs $G^{n,L}$ are both irreducible and recurrent,
so that the stationary distribution is unique (up to normalization).
We will use $M^{n,L}$ to denote the transition matrix of our Markov chain,
	\be
	M^{n,L}_{v,u} = \text{rate}(u \to v) = \beta(v).
	\ee
$\wM^{n,L}$ is then precisely the transition matrix,
	\be
	{}^\nabla\! M^{n,L}_{v,u} = \begin{cases}
	\beta(v) & \text{for}~u \neq v, \\
	 \ds -\sum_{\substack{w \in \Sigma^{L} \\ u \neq w}}\beta(w) & \text{for}~u=v.
	\end{cases}
	\ee
For the example in Figure~\ref{fig:example}, with lexicographic ordering of the states,
	\be
	{}^\nabla\!M^{2,2}=  
	\left( \begin {array}{cccc} 
	-x_{2,{1}}&0&x_{1,{2}}&0\\ 
	\noalign{\medskip}
	x_{2,{1}}&-x_{1,{1}}-x_{2,{2}}&x_{2,{1}}&0\\ 
	\noalign{\medskip}
	0&x_{1,{1}}&-x_{1,{2}}-x_{2,{1}}&x_{1,{1}}\\ 
	\noalign{\medskip}
	0&x_{2,{2}}&0&-x_{1,{1}}
	\end {array} \right).
	\ee
The stationary distribution is given by probabilities of words, 
which are to be taken as rational functions in the variables $x_{a,t}$.
It is  the column vector with eigenvalue zero, which after normalization is then given by
	\be
	\begin{split}
	&\text{Pr}[1,1]={\frac {x_{1,{1}}x_{1,{2}}}
	{ \left( x_{1,{2}}+x_{2,{1}} \right)  \left( x_{1,{1}}+x_{2,{1}} \right) }},
	\, \text{Pr}[1,2]={\frac {x_{2,{1}}x_{1,{1}}}
	{ \left( x_{1,{1}}+x_{2,{2}} \right)  \left( x_{1,{1}}+x_{2,{1}} \right) }},\\
	&\text{Pr}[2,1]={\frac {x_{2,{1}}x_{1,{1}}}
	{ \left( x_{1,{2}}+x_{2,{1}} \right)  \left( x_{1,{1}}+x_{2,{1}} \right) }},
	\, \text{Pr}[2,2]={\frac {x_{2,{2}}x_{2,{1}}}
	{ \left( x_{1,{1}}+x_{2,{2}} \right)  \left( x_{1,{1}}+x_{2,{1}} \right) }}.  
	\end{split}
	\ee
Notice that the probabilities consist of a product of two monomials in the numerator
and two factors in the denominator, and that each factor contains two terms. Also, 
notice that not all the denominators are the same, otherwise the steady state would 
be a true product measure. Of course, the sums of these probabilities is 1, which is 
not completely obvious.

It is also interesting to note that the eigenvalues of ${}^\nabla\!M^{2,2}$ are linear
in the variables. Other than zero, the eigenvalues are given by
	\be
	-x_{1,{1}}-x_{2,{2}}, \quad -x_{1,{1}}-x_{2,{1}}, \text{ and } -x_{1,{2}}-x_{2,{1}}.
	\ee
Another way of saying this is that the characteristic polynomial of the 
transition matrix factorizes into linear parts.

\section{Stationary Distribution}
In this section we determine an explicit expression for the steady state distribution
of the de Bruijn process on $G^{n,L}$. Before we do that we will have to set down 
some notation. For the moment we are working over $\Sigma^{+}$, the set of all 
nonempty words over the alphabet $\Sigma$ (of size $n$). 

For convenience, we introduce operators which denote the transitions of our Markov chain.
Let $\partial_{a}$ be the operator that adds the letter $a$ to the end of a word and 
removes the first letter,
	\be
	\partial_{a} : u \mapsto  t(u). a.
	\ee
With $\beta$ as introduced we introduce the shorthand notation
	\be
	\beta_{a,m} =  \sum_{ b \in \Sigma} \beta(\partial_b \,a^m) = 
		x_{a,m} + \sum_{ b \in \Sigma, b \neq a}Êx_{b,1}.
	\ee
Note that $\beta_{a,1}= \sum_{b \in \Sigma} x_{b,1}$ does not depend on $a$.
We now define the valuation $\mu(u)$ for $u \in \Sigma^+$ as
	\be \label{defmu}
	\mu(u) = \frac{\beta(u)}{ \sum_{a \in \Sigma} \beta(\partial_{a}u)}.
	\ee
Note that the restriction of $\mu$ on the alphabet $\Sigma$ is (formally) a probability distribution.
Finally, we define the valuation $\bar\mu$, also on $\Sigma^{+}$, as
	\be \label{defbarmu}
	\bar\mu(u) = \prod_{i=1}^{L} \mu(u_{1} u_{2} \dots u_{i}) 
	= \mu(u_{1}) \mu(u_{1} u_{2}) \cdots \mu(u_{1} u_{2} \dots u_{L}),
	\ee
if $u=u_1u_2\ldots u_L$.	
The following result is the key to understanding the stationary distribution.
\begin{prp} \label{prp:barmu}
For all $u \in \Sigma^{+}$,
	\be \label{barmuiden}
	\sum_{ a \in \Sigma} \bar\mu(a . u) = \bar\mu(u).
	\ee
\end{prp}

\begin{proof}
As in (\ref{blockfac}), let us write $w$ in block factorized form:
	\be\label{blockfactor}
	u=b^{(1)} b^{(2)} \cdots b^{(m)} = \tilde w . b^{(m)},
	\ee
where  $\tilde{u}=b^{(1)}Ê\ldots b^{(m-1)}$ if $m>1$, and 
$\tilde{u}$ is the empty word if $m=1$.

If $b^{(m)}=a^{k}$,
then
	\be
	\mu(u) = \begin{cases}
	 \ds \frac{x_{a,k}}{ \beta_{a,k}} & \text{if}~m=1, \text{i.e., if $u$ is a block}, \\[5mm]
	 \ds \frac{x_{a,k}}{\beta_{a,k+1}} & \text{if}~m>1,
	\end{cases}
	\ee
and thus
	\be
	\bar\mu(u) = \begin{cases}
	\ds \prod_{j=1}^{k} \frac{x_{a,j}}{ \beta_{a,j}} & \text{if}~m=1, \text{i.e., if $u$ is a block}, \\
	\ds \bar\mu(\tilde w) \cdot \prod_{j=1}^{k} \frac{x_{a,j}}{ \beta_{a,j+1}} & \text{if}~m>1.
	\end{cases}
	\ee

We will define another valuation on $\Sigma^{+}$ closely related to $\bar\mu$, 
which we call $\bar\rho$.
Referring to the factorization (\ref{blockfactor}) we put
	\be
	\bar\rho(u) = \begin{cases}
	\ds \prod_{j=1}^{k} \frac{x_{a,j}}{\beta_{a,j+1}} & 
		\text{if}~m=1, \text{i.e., if $u=a^k$ is a block}, \\
	\ds  \prod_{l=1}^{m} \bar\rho(u^{(l)}) & \text{if}~m>1.
	\end{cases}
	\ee
This new valuation is related to $\bar\mu$ by the following properties:
\begin{itemize}
\item[--] For blocks $u=a^k$ we have
	\be \label{prop1}
	\bar\rho(a^{k}) = \frac{ \beta_{a,1}}{ \beta_{a,k+1}} \bar\mu(a^{k}),
	\ee
\item[--] For $u$ with factorization (\ref{blockfactor}) we have
	\be \label{prop2}
	\bar\mu(u) = \bar\mu(\tilde u) \cdot  \bar\rho(b^{(m)}),
	\ee
\item[--]
which, by the obvious induction, implies
	\be
	\bar\mu(u) = \bar\mu(b^{(1)}) \cdot  \prod_{l=2}^{m} \bar\rho(b^{(l)}).
	\ee
\end{itemize}

We are now in a position to prove identity \eqref{barmuiden}.
First consider the case where $u=a^k$ is a block.
	
	\be
	\begin{split}
	\sum_{b \in \Sigma} \bar\mu(b \cdot a^{k}) &= \bar\mu(a^{k+1}) 
	+ \sum_{b \neq a} \bar\mu(b \cdot a^{k}) \\
	&=\frac{x_{a,k+1}}{ \beta_{a,k+1}}  \bar\mu(a^{k}) 
	+ \sum_{b \neq a} \bar\mu(b) \cdot \bar\rho( a^{k}) \\
	&=\frac{x_{a,k+1}}{ \beta_{a,k+1}}  \bar\mu(a^{k}) 
	+ \sum_{b \neq a} \frac{x_{b,1}}{\beta_{a,1}} \bar\rho(a^{k}) \\
	&=\left( \frac{x_{a,k+1}}{ \beta_{a,k+1}} + 
	\sum_{b \neq a} \frac{x_{b,1}}{ \beta_{a,k+1}} \right)
	  \bar\mu(a^{k}) \\
	  &= \bar\mu(a^{k}),
	\end{split}
	\ee
where we used \eqref{prop1} in the last-but-one step.
The general case is then proven by a simple induction on $m$.
	\be
	\begin{split}
	\sum_{a \in \Sigma} \bar\mu(a . b^{(1)} b^{(2)} \ldots b^{(m)}) &=
	\sum_{a \in \Sigma} \bar\mu(a . b^{(1)} b^{(2)} \ldots b^{(m-1)}) \cdot \bar\rho(b^{(m)}) \\
	&= \bar\mu(b^{(1)} b^{(2)} \ldots b^{(m-1)}) \cdot \bar\rho(b^{(m)}) \\
	&= \bar\mu(b^{(1)} b^{(2)} \ldots b^{(m)}),
	\end{split}
	\ee
where we have used  property \eqref{prop2} of $\bar\rho$ in the last step.
\end{proof}

As a consequence of Proposition~\ref{prp:barmu}, we have the following result, which is
an easy exercise in induction. 
The case $L=1$ was already mentioned immediately after \eqref{defmu}.

\begin{cor}
For any fixed length $L$ of words over the alphabet $\Sigma$,
	\be
	\sum_{w \in \Sigma^{L}} \bar\mu(w) = 1.
	\ee
\end{cor}
Therefore, the column vector 
$\boldsymbol{\bar\mu}^{n,L} = [\bar\mu(u)]_{u \in \Sigma^{L}}$ 
can be a seen as a formal
probability distribution on $\Sigma^{L}$.
We now look at the transition matrix $M^{n,L}$ more closely. 
	\be
	M^{n,L}_{v,u} = \delta_{h(v)=t(u)} \, \beta(v).
	\ee
where $\delta_{x}$ is the indicator function for $x$, i.e.,
it is 1 if the statement $x$ is true and 0 otherwise.
Thus the matrix $M^{n,L}$ is very sparse. It has just $n$ off-diagonal 
non-zero entries per row and per column. 
More precisely, the row indexed by $v$ has the entry $\beta(v)$ for 
the $n$ $\partial$-preimages of $v$, and the column indexed by 
$w$ contains $\beta(\partial_{a} u)$ as the only nonzero entries. 
In particular, the  column sum for the column indexed by $u$ is 
$\sum_{a \in \Sigma} \beta(\partial_{a}(u))$.
Define the diagonal matrix $\Delta^{n,L}$ as one with with precisely 
these column sums as entries, i.e.
	 \be
	 \Delta^{n,L}_{v,u} = \begin{cases}
	 \sum_{a \in \Sigma} \beta(\partial_{a}u) & v=u, \\
	 0 & \text{otherwise}.
	 \end{cases}
	 \ee

\begin{thm} \label{thm:stat}
The vector $\bmubar^{n,L}$ is the stationary vector for the de Bruijn 
process on $G^{n,L}$, i.e.,
	\be \label{stat}
	M^{n,L} \boldsymbol{\bar\mu}^{n,L} = \Delta^{n,L} \boldsymbol{\bar\mu}^{n,L}.
	\ee
\end{thm}

\begin{proof}

Consider the row corresponding to word 
$v = v_1v_2 \ldots v_{L-1}v_L= h(v). v_L$ in
the equation
	\be\label{stationary}
	M \, \boldsymbol{\overline{\mu}} = \Delta \, \boldsymbol{\overline{\mu}}.
	\ee
On the l.h.s. of (\ref{stationary}) we have to consider
the summation  
$\sum_{u \in \Sigma^L} M_{v,u}\, \overline{\mu}(u)$,
where only those $u u \in \Sigma^L$ with $t(u).v_L=v$
contribute. 
This latter condition can be written as 
$u=b.h(v)$ for some $b \in \Sigma$,
so that this summation can be written as
	\be
	\begin{split}
	\ds \sum_{u \in \Sigma^L} M_{v,u}\, \overline{\mu}(u) b&= 
	\ds \sum_{b \in \Sigma} M_{v,b . h(v)} \, \overline{\mu}(b . h(v)) \cr
	&= \beta(v) \ds \sum_{b \in \Sigma} \overline{\mu}(b . h(v)) =
	\beta(v) \, \overline{\mu}(h(v)),
	\end{split}
	\ee
where the last equality follows from Lemma 1.
	
On the r.h.s. of  (\ref{stationary}) 	
we have for the row entry corresponding to the word $v$:
	\be
	\begin{split}
	\Delta_{v,v}\, \overline{\mu} (v) &= \ds \sum_{a \in \Sigma} \beta(\partial_a v) \, \overline{\mu}(v) \cr
	&= \ds \sum_{a \in \Sigma} \beta(\partial_a v) \cdot  \overline{\mu}(h(v)) \,\mu(v) 
	=\beta(v) \, \overline{\mu}(h(v))
	\end{split}
	\ee
in view of the inductive definition of $\overline{\mu}$ in (\ref{defbarmu})
and the definition of $\mu$ in (\ref{defmu}). 
\end{proof}

Let $Z^{n,L}$ denote the common denominator of  the stationary probabilities of
configurations. This is often called, with some abuse of terminology, 
the {\em partition function} \cite{blytheevans}. The abuse comes from the fact 
that this terminology is strictly applicable in the sense of statistical mechanics 
while considering Markov chains only when they are reversible. 
The de Bruijn process definitely does not fall into this category.
Since the probabilities are given by products of $\mu$ in \eqref{defbarmu}, 
one arrives at the following product formula.
\begin{cor} \label{cor:partfn}
The partition function of the de Bruijn process on $G^{n,L}$
is given by	
	\be
	Z^{n,L} = \beta_{1,1} \cdot  \prod_{m=2}^{L-1} \prod_{a=1}^{n} 
	 \beta_{a,m}.
	\ee
\end{cor}

Physicists are often interested in properties of the stationary distribution 
rather than the full distribution itself. One natural quantity of interest in this 
context is the so-called density distribution of a particular letter, say $a$, in the alphabet. 
In other words, they would like to know, for example, how likely it is that 
$a$ is present at the first site rather than the last site. We can make this 
precise by defining {\em occupation variables}. 
Let $\eta^{a,i}$ denote the occupation variable of species $a$ at site $i$: 
it is a random variable which is 1 when site $i$ is occupied by $a$ and zero otherwise.
We define the probability in the stationary distribution by the symbol 
$\langle \;\cdot \; \rangle$.
Then $\langle \;\eta^{a,i} \; \rangle$ gives the {\em density} of $a$ at site $i$.
Similarly, one can ask for joint distributions, such as 
$\langle \;\eta^{a,i} \eta^{b,j} \; \rangle$, which is the probability that site $i$ 
is occupied by $a$ and simultaneously that site $j$ is occupied by $b$. 
Such joint distributions are known as {\em correlation functions}.

We will not be able to obtain detailed information about arbitrary correlation 
functions in full generality, but there is one case in which we can easily give the answer.
This is the correlation function for any letters $a_{k},  \dots, a_{2}, a_{1}$ at the last $k$ sites.
\begin{cor} \label{cor:occuplastk}
Let $u=a_{k}  \ldots a_{2} a_{1}$. Then
	\be
	\langle \eta^{a_{k},L-k+1} \cdots  \eta^{a_{2},L-1} \eta^{a_{1},L} \rangle
	= \bar \mu(u).
	\ee
\end{cor}

\begin{proof}
By definition of the stationary state,
	\be
	\langle \eta^{a_{k},L-k+1} \cdots  \eta^{a_{2},L-1} \eta^{a_{1},L} \rangle
	= \ds \sum_{v \in \Sigma^{L-k}} \bar\mu(v . u).
	\ee
Using Proposition~\ref{prp:barmu} repeatedly $L-k$ times, we arrive at the desired result.
\end{proof}

In particular, Corollary~\ref{cor:occuplastk} says that the density of species $a$ at 
the last site is simply
	\be
	\langle \eta^{a,L} \rangle = \frac{ x_{a,1} }{ \beta_{a,1}}.
	\ee
Formulas for densities at other locations are much more complicated. 
It would be interesting to find a uniform formula for the density of species $a$ at site $k$.

\section{Characteristic Polynomial of $\wM^{n,L}$}
We will prove a formula for the characteristic polynomial of $\wM^{n,L}$ in the following. 
In particular, we will show that it factorizes completely into linear parts.
In order to do so, we need to understand the structure of the transition matrices better.
We denote by $\chi(M;\lambda)$ the characteristic polynomial of a matrix $M$ 
in the variable $\lambda$. 

To begin with, let us recall from the previous section
that the transition matrices $M^{n,L}$, taken as
mappings defined on row and column indices, are defined by
	\be
	M^{n,L}:\Sigma_n^L \times \Sigma_n^L \rightarrow X : 
	(v,u) \mapsto \delta_{h(v)=t(u)} \cdot \beta(v).
	\ee

\begin{lem}
The matrix $M^{n,L}$ can be written as
	\be
	M^{n,L} = \left[ \, A^{n,L} \,|\, A^{n,L} \,|\, \ldots \,|\, A^{n,L} \, \right]
	~~~~(n~\text{copies of}~A^{n,L}),
	\ee
where $A^{n,L}$ is a matrix of size $n^L \times n^{L-1}$ given by
	\be
	A^{n,L} : \Sigma^{n,L} \times \Sigma^{n,L-1} \rightarrow X \cup \{0\} : 
	(v,u) \mapsto \delta_{h(v)=u} \cdot \beta(v).
	\ee
\end{lem}

We have
	\be \label{recura}
	 A^{n,1} = 
	 \begin{bmatrix}
	 x_{1,1} \cr x_{2,1} \cr \vdots \cr x_{n,1}
	 \end{bmatrix},\hspace{0.5cm}
	A^{n,L} = 
	 \begin{bmatrix}
	A^{n,L-1}_{1} & 0^{n,L-1} & \cdots  & 0^{n,L-1} \\
	0^{n,L-1} & A^{n,L-1}_{2} & \cdots  & 0^{n,L-1} \\
	\vdots & \vdots & \ddots & \vdots\\
	0^{n,L-1} & 0^{n,L-1} & \cdots  & A^{n,L-1}_{n}
	 \end{bmatrix}=
	 \begin{bmatrix} B_1^{n,L-1} \cr B_2^{n,L-1} \cr \vdots \cr B_{n}^{n,L-1}
				\end{bmatrix},
	 \ee
where $A^{n,L-1}_{k}$ is like $A^{n,L-1}$, but with $x_{k,L-1}$ replaced by $x_{k,L}$, 
and where $0^{n,L-1}$ is the  zero matrix of size $n^{L-1} \times n^{L-2}$.
The matrices $B_a^{n,L-1}$ are square matrices of size $n^{L-1} \times n^{L-1}$, where
for each $a\in \Sigma$  the matrix $B_a^{n,L}$ is defined by
	\be
	B_a^{n,L} : \Sigma^L \times \Sigma^{L} \rightarrow X \cup \{0\} :
	(v,u) \mapsto \delta_{a.h(v)=u} \cdot \beta(a.v).
	\ee
With these matrices at hand we can finally define the matrix 
$B^{n,L}=\sum_{a \in \Sigma} B_a^{n,L}$ of size $n^{L} \times n^{L}$, so that
	\be
	B^{n,L} : 
	\Sigma^L \times \Sigma^{L} \rightarrow X\cup \{0\} :
	(v,u) \mapsto \delta_{h(v)=t(u)} \cdot \beta(u_1.v).
	\ee

\begin{lem}
	$M^{n,L}- B^{n,L}$ is a diagonal matrix.
\end{lem}

\begin{proof}	We have
	\be
	M^{n,L}(v,u) \neq B^{n,L}(v,u) ~~\Leftrightarrow~~
	h(v)=t(u)~\text{and}~\beta(u_1.v)\neq \beta(v)
	\ee
But $\beta(u_1.v)\neq \beta(v)$ can only happen if the 
last block of $u_1.v$ is different from the last block of $v$,
which only happens if $v$ itself is a block, $v=a^L$,
and $u_1=a$, in which case $\beta(v)=x_{a,L}$
and $\beta(u_1.v)=x_{a,L+1}$. 
So we have
	\be
	(B^{n,L}-M^{n,L})(v,u) = 
	\begin{cases}
	x_{a,L+1}-x_{a,L} & \text{if}~v=u=a^L ,\cr
	0 & \text{otherwise.}
	\end{cases}
	\ee
\end{proof}
We state as an equivalent assertion:

\begin{cor}\label{cor}
For the Kirchhoff matrices of $M^{n,L}$ and $B^{n,L}$ we have equality: 
	\be
	{}^\nabla\!M^{n,L} ={}^\nabla\!B^{n,L}.
	\ee
\end{cor}

We now prove a very general result about the characteristic polynomial 
of a matrix with a certain kind of block structure. 
This will be the key to finding the characteristic polynomial of our transition matrices.

\begin{lem} \label{lem:blockm}
Let $P_{1},\dots,P_{m},Q$ be any $k \times k$ matrices, $P = P_{1}+ \dots+P_{m}$ and 
	\be
	R= \begin{bmatrix} 
	P_{1}+Q & P_{2} & \cdots  & P_{m} \\
	P_{1} & P_{2}+Q & \cdots  & P_{m} \\
	\vdots & \vdots & \ddots & \vdots \\
	P_{1} & P_{2} & \cdots  & P_{m}+Q \\
	 \end{bmatrix}.
	\ee
Then
	\be
	\chi(R;\lambda) =  \chi(Q;\lambda)^{m-1} \cdot \chi(P + Q;\lambda).
	\ee
\end{lem}

\begin{proof}
Multiply $R$ by the block lower-triangular matrix of unit determinant shown to get
	\be
	R \cdot \begin{bmatrix}
	1 & 0 & 0 & \cdots  & 0 \\
	-1 & 1 & 0 & \cdots  & 0 \\
	0  & -1 & 1 & \cdots  & 0 \\
	\vdots & \vdots & \vdots & \ddots & \\
	0 & 0 & 0  & \cdots & 1 \\
	 \end{bmatrix}
	 =
	 \begin{bmatrix} 
	Q & 0 & 0 & \cdots & P_{1} \\
	-Q & Q & 0 & \cdots  & P_{2} \\
	0  & -Q & Q & \cdots & P_{3} \\
	\vdots & \vdots & \vdots &\ddots & \vdots \\
	0 & 0 & 0  & \cdots & P_{m}+Q \\
	 \end{bmatrix}
	 \ee
which has the same determinant as $R$.
Now perform the block row operations which replace row $j$ by the sum of
rows 1 through $j$ to get
	\be
	 \begin{bmatrix} 
	Q & 0 & 0 & \cdots & P_{1} \\
	0 & Q & 0 & \cdots  & P_{1}+P_{2} \\
	0  & 0 & Q & \cdots & P_{1}+P_{1}+P_{3} \\
	\vdots & \vdots & \vdots &\ddots & \vdots \\
	0 & 0 & 0  & \cdots & P+Q \\
	 \end{bmatrix}
	 \ee
Since this is now a block upper triangular matrix, the characteristic 
polynomials is the product of those of the diagonal blocks.
\end{proof}

We will now apply this lemma to the block matrix
	\be
	\wM^{n,L+1} =  
	\begin{bmatrix}
	B_1^{n,L} - D^{n,L} & B_1^{n,L} & \ldots & B_1^{n,L} \cr
	B_2^{n,L} & B_2^{n,L} -D^{n,L} & \ldots & B_2^{n,L} \cr
	\vdots & \vdots & \ddots & \vdots \cr
	B_{n}^{n,L} & B_{n}^{n,L}& \ldots & B_{n}^{n,L} -D^{n,L}
	\end{bmatrix}
	\ee
where $D^{n,L}$ is the $(n^{L} \times n^{L})$-diagonal matrix
	with the column sums of $A^{n,L+1}$ on the main diagonal. 

\begin{prp}\label{prop}
	The characteristic polynomials $\chi(\wM^{n,L};z)$ satisfy the recursion
	\be
	\chi(\wM^{n,L+1};z) =\chi(-D^{n,L};z)^{n-1} \cdot \chi(\wM^{n,L};z).
	\ee
\end{prp}

\begin{proof}
From Corollary \ref{cor}, Lemma \ref{lem:blockm},  and the easily checked fact 
${}^\nabla\!B^{n,L} = B^{n,L} - D^{n,L}$ we get:
	\be
	\begin{split}
	\chi(\wM^{n,L+1};\lambda) &= 
	\chi(-D^{n,L};\lambda)^{n-1} \cdot \chi(\ts \sum_{a \in \Sigma} B_a^{n,L}-D^{n,L};\lambda) \cr
	&=
	\chi(-D^{n,L};\lambda)^{n-1} \cdot \chi(B^{n,L}-D^{n,L};\lambda) \cr
	&=
	\chi(-D^{n,L};\lambda)^{n-1} \cdot \chi(\wB^{n,L};\lambda) \cr
	&=
	\chi(-D^{n,L};\lambda)^{n-1} \cdot \chi(\wM^{n,L};\lambda).
	\end{split}
	\ee
\end{proof}

 As a final step, we need a formula for $\chi( - D^{n,L}, \lambda)$.
\begin{lem} \label{lem:detx}
The characteristic polynomial of $-D^{n,L}$ is given by
	\be
	\chi( - D^{n,L}, \lambda) = \begin{cases}
	\ds \lambda + \beta_{1,1}  &  \text{if}~L=1, \\
	\\
	 \ds \prod_{m=2}^{L-1} \prod_{a \in \Sigma} 
	 \left(  \lambda +\beta_{a,m} \right)^{(n-1) n^{L-1-m}} 
	\; \ds \prod_{a \in \Sigma} 
	\left( \lambda + \beta_{a,L} \right) 
	& \text{if}~L>1.
	\end{cases}
	\ee
\end{lem}

\begin{proof}
The case $L=1$ follows directly from the definition of $A^{n,1}$ in \eqref{recura}.
For general $L$, recall that $A^{n,L}$ contains $n$ copies of $A^{n-1,L}$
with one factor containing $x_{a,L-1}$ removed and one factor containing $x_{a,L}$ 
added instead, for each $a \in \Sigma$.
Thus,
	\be
	\chi(-D^{n,L},\lambda) = \left[ \chi(-D^{n,L-1},\lambda) \right]^{n} \cdot
	\prod_{a \in \Sigma} \left( \frac{\lambda+ \beta_{a,L} }
	{\lambda+\beta_{a,L-1} } \right),
	\ee
which proves the result.
\end{proof}

We can now put everything together and get from Proposition \ref{prop}, Lemma \ref{lem:detx}
and checking the initial case for $L=1$:
\begin{thm} \label{thm:charpoly}
The characteristic polynomial of the de Bruijn process on $G^{n,K}$ is given by
	\be \label{detm}
	\chi(\wM^{n,L}; \lambda) =\lambda \,(\lambda + \beta_{1,1}) \cdot 
	\prod_{m=2}^{L} \prod_{a \in \Sigma} 
	 \left(  \lambda + \beta_{a,m} \right)^{(n-1)n^{L-m}}.
	\ee
\end{thm}

\section{Special cases}
We now consider special cases of the rates where something 
interesting happens in the de Bruijn process.

\subsection{The de Bruijn-Bernoulli Process}
There turns out to be a special case of the rates $x_{a,j}$ for 
which the stationary distribution is a {\em Bernoulli measure}. 
That is to say, the probability of finding species $a$ at site 
$i$ in stationarity is independent, not only of any other site, but also of $i$ itself.
This is not  obvious because the dynamics at any given site is certainly a priori not
independent from what happens at any other site. Since the measure is so simple, 
all correlation functions are trivial. We denote the single site measure in \eqref{defmu} 
for this specialized process to be $\mu_{y}$, and the stationary measure \eqref{defbarmu} 
as $\bar\mu_{y}$.

\begin{cor} \label{cor:prod}
Under the choice of rates $x_{a,j} = y_{a}$ independent of $j$, 
the stationary distribution of the Markov chain with transition 
matrix $\wM^{n,L}$ is Bernoulli with density
	\be \label{defrho}
	\rho_{a} = \frac{y_{a}}{\ts\sum_{b \in \Sigma} y_{b}}.
	\ee
\end{cor} 

\begin{proof}
The choice of rates simply mean that species $a$ is added 
with a rate independent of the current configuration. From \eqref{defmu},
it follows that for $u=u_1u_2\ldots u_L$,
	\be
	\mu_{y}(u) = \frac{y_{u_{L}}}{\ts\sum_{b \in \Sigma} y_{b}} = \rho_{u_{L}},
	\ee
and using the definition of the stationary distribution $\bar\mu$ in \eqref{defbarmu},
	\be
	\bar\mu_{y}(u) = \prod_{i=1}^{L } \rho_{u_{i}},
	\ee
which is exactly the definition of a Bernoulli distribution.
\end{proof}

\subsection{The Skin-deep de Bruijn Process}
 
Another tractable version of the de Bruijn process is one where 
the rate for transforming the word $u=u_1u_2\ldots u_L$
into  $\partial_a u= t(u).a=u_2\ldots u_L.a$ for $a \in \Sigma$
only depends on the occupation of the last site,  $u_{L}$. 
Hence, the rates are only {\em skin-deep}.
An additional simplification comes by choosing the rate to be 
$x$ when $a=u_{L}$ and 1 otherwise. Namely,
	\be
	x_{a,j} = \begin{cases}
	x & \text{for}~j=1, \\
	1 & \text{for}~j>1.
	\end{cases}
	\ee

We first summarize the results. It turns out that any letter in the alphabet 
is equally likely to be at any site in the skin-deep de Bruijn process. 
This is an enormous simplification compared to the original process where 
we do not have a general formula for the density. Further, we have the 
property that all correlation functions are independent of the length of the words. 
This is not obvious because the Markov chain on words of length $L$ 
is not reducible in any obvious way to the one on words of length $L-1$.
This property is quite rare and very few examples are known of such families of Markov chains. 
One such example is the asymmetric annihilation process \cite{ayyerstrehl}.

The intuition is as follows. By choosing $x \ll 1$ one prefers to add 
the same letter as $u_{L}$, and similarly, for $x \gg 1$, one prefers 
to add any letter in $\Sigma$ other than $u_{L}$. 
Of course, $x=1$ corresponds to the uniform distribution. 
Therefore, one expects the average word to be qualitatively different 
in these two cases. {\em In the former case, one expects the 
average word to be the same letter repeated $L$ times, 
whereas in the latter case, one would expect  no two neighboring letters 
to be the same on average.}
Our final result, a simple formula for the two-point correlation function, 
exemplifies the different in these two cases.

We begin with a formula for the stationary distribution, which we 
will denote in this specialization by $\bar\mu_{x}$.
We will always  work with the alphabet $\Sigma$ on $n$ letters.

\begin{lem} \label{lem:spbarmu}
The stationary probability for a word $u=u_1u_2\ldots u_L \in \Sigma^L$  is 
given by 
	\be \label{spbarmu}
	\bar\mu_x(u) = \frac{x^{ \gamma(u)-1}}
	{ n (1+(n-1)x)^{L-1}},
	\ee
where $\gamma(u)$ is the number of blocks of $u$.
\end{lem}

\begin{proof}
Analogous to the notation for the stationary distribution, 
we denote the block function by $\beta_{x}$.
From the definition of the model, 
	\be \label{spbeta}
	\beta_x(a^{k}) = \begin{cases}
	x & \text{if}~k=1,\\
	1 & \text{if}~k>1.
	\end{cases}
	\ee
and thus, for any word $u$ the value 
$\beta_x(u)$ is $x$ if the length of the last block in its block
decomposition is 1, and is $1$ otherwise. 
The denominator in  \eqref{spbarmu} is easily explained. 
For any word $u$ of length $L$,
	\be
	\ds \sum_{a \in \Sigma}\beta_x(t(u) . a) = \begin{cases}
	1+(n-1)x & L>1,\\
	nx & L=1,
	\end{cases}
	\ee
because for all but one letter in $\Sigma$, the size of the last 
block in $t(u).a$ is going to be 1. 
The only exception to this argument is, $L=1$, when $t(u)$ is empty. 
From \eqref{defbarmu}, we get
	\be
	\bar\mu_x(u) = \frac{ \beta_x(u_{1}) \beta_x(u_{1} u_{2}) \cdots \beta_x(u_{1}\dots u_{L})}
	{ nx (1+(n-1)x)^{L-1}}.
	\ee
The numerator is $x^{\gamma(u)}$, since we pick up a factor of $x$ every time 
a new block starts.  One factor $x$ is cancelled because $\beta_x(u_{1})=x$.
\end{proof}

The formula for the density is essentially an argument about the symmetry of the de Bruijn graph $G^{n,L}$.

\begin{cor} \label{cor:spdens}
The probability in the stationary state of $G^{n,L}$ that site $i$ is occupied by letter $a$ is uniform,
i.e., for any $i$ s.th. $1 \leq i \leq L$ we have 
	\be
	\langle \eta^{a,i} \rangle = \frac 1n ~~~(a \in \Sigma).
	\ee
\end{cor}

\begin{proof}
Indeed, by Lemma~\ref{lem:spbarmu} the stationary distribution $\bar{\mu}_x$
is invariant under any permutation of the letters of the alphabet $\Sigma$.
Hence  $\langle \eta^{a,i} \rangle$ does not depend on $a \in \Sigma$ and
we have uniformity.
\end{proof}

Since the de Bruijn-Bernoulli process has a product measure, 
the density of $a$ at site $i$ is also independent of $i$, 
but the density is not uniform since it is given by $\rho_{a}$
\eqref{defrho}. The behavior of higher correlation functions 
here is more complicated than the de Bruijn-Bernoulli process. 
There is, however, one  aspect in which it resembles the former,
namely:

\begin{lem} \label{lem:spindofL}
Correlation functions of $G^{n,L}$ in this model
are independent of the length $L$ of the words and they are shift-invariant. 
\end{lem}

\begin{proof}
We can represent an arbitrary correlation function in the de Bruijn graph $G^{n,L}$ as 
	\be\label{eq:correlation}
	\langle \eta^{a_{1},i_{1}} \cdots \eta^{a_{k},i_{k}} \rangle_L = 
	\sum_{w^{(0)}, \dots, w^{(k)}} \bar\mu_x(w^{(0)}a_{1}w^{(1)} \dots w^{(k-1)} a_{k} w^{(k)}),
	\ee
where we have sites 
$1 \leq i_1 < i_2 < \ldots < i_k \leq L$ and letters $a_1,a_2,\ldots,a_k \in \Sigma$, 
and where the sum runs over all $(w^{(0)},w^{(1)},\ldots,w^{(k)})$ with
$w^{(j)} \in \Sigma^{i_{s+1}-i_{s}-1}$ for $s \in \{0,\dots,k \}$,
and where we put $i_0=0$ and $i_{k+1}=L+1$.
Now note that we have from Proposition \ref{prp:barmu} for any $u \in \Sigma^k$
	\be
	\sum_{w \in \Sigma^\ell}Ê\bar\mu_x (w.u) = \bar\mu_x(u).
	\ee
Since $\bar\mu_x$, as given in Lemma \ref{lem:spbarmu}, is also invariant under reversal
of words, we also have
	\be
	\sum_{w \in \Sigma^\ell}Ê\bar\mu_x (u.w) = \bar\mu_x(u).
	\ee
As a consequence, we can forget about the outermost summations in
(\ref {eq:correlation}) and get
	\begin{multline}\label{eq:correlation}
	\langle \eta^{a_{1},i_{1}} \cdots \eta^{a_{k},i_{k}} \rangle_L = \cr
	\sum_{w^{(1)}, \dots, w^{(k-1)}} \bar\mu_x(a_{1}w^{(1)} \dots w^{(k-1)} a_{k}) =
	\langle \eta^{a_{1},j_{1}} \cdots \eta^{a_{k},j_{k}} \rangle_{i_k-i_1+1},
	\end{multline}
where $j_s=i_s-i_1+1~(1 \leq s \leq k)$.
	Shift-invariance in the sense that
	\be
	\langle \eta^{a_{1},i_{1}} \cdots \eta^{a_{k},i_{k}} \rangle_L
	=\langle \eta^{a_{1},i_{1}+1} \cdots \eta^{a_{k},i_{k}+1} \rangle_L
	\ee
is an immediate consequence.
\end{proof}

We now proceed to compute the two-point correlation function. This is an
easy exercise in generating functions for words according to the number of blocks.
The technique is known as ``transfer-matrix method'', see, e.g., Section 4.7
in \cite{stanley1}.

For $a,b \in \Sigma$ and $k \geq 1$ we define the generating polynomial in the variable $x$
	\be\label{eq:defalpha}
	\alpha_{n,k}(a,b;x)=\sum_{w \in a.\Sigma^{k-1}.b}   x^{\gamma(w)-1} ,
	\ee
where, as before, $\gamma(w)$ denotes the number of blocks in the 
block factorization of $w\in \Sigma^+$ (so that $\gamma(w)-1$ is the 
number of pairs of adjacent distinct letters in $w$).
Note that
	\be
	\alpha_{n,1}(a,b;x) = \begin{cases} 1 &\text{if}~a=b, \cr 
								x & \text{if}~a \neq b. \end{cases}
	\ee
 The following statement is folklore:
\begin{lem}\label{lem:folklore}
	Let $\mathbb{I}_n$ denote the identity matrix
        and $\mathbb{J}_n$ denote the all-one matrix, both of size $n \times n$, and let
	$K_n(s,t) := s \cdot \mathbb{I}_n + t\cdot \mathbb{J}_n$ for parameters $s,t$. 
	Then 
	\be
	K_n(s,t)^{-1} = \frac{1}{s(s+nt)} \, K_n(s+nt,-t).
	\ee
\end{lem}
Indeed, this is a very special case of what is known as the Sherman-Morrison formula,
see \cite{sm50}, \cite{wilf59}.

 Consider now the matrix 
	\be
	A_n(x) := \left[\, \alpha_{n,1}(a,b;x) \,\right]_{a,b \in \Sigma} =
	(1-x)\cdot \mathbb{I}_n + x \cdot \mathbb{J}_n = K_n(1-x,x)
	\ee
which encodes transition in the alphabet $\Sigma$. Then, for $k \geq 1$,
$A_n(x)^k$ is an $(n \times n)$-matrix which in position $(a,b)$
contains the generating polynomial $\alpha_{n,k}(a,b;x)$:
	\be
	A_n(x)^k = \left[ \, \alpha_{n,k}(a,b;x) \, \right]_{a.b \in \Sigma}.
	\ee
We can get generating functions by summing the geometric series
and using Lemma \ref{lem:folklore}:
	\be
	\begin{split}
	\ds\sum_{k \geq 0} A_n(x)^k z^k&= (\mathbb{I}_n - z \cdot A_n(x))^{-1}  \cr
	&= K_n(1-z+xz,-xz)^{-1} \cr
	&= \frac{K_n(1-z-(n-1)xz,xz)}{(1-z+xz)(1-z-(n-1)xz)}, 
	\end{split}
	\ee
which means that that for any two distinct letters $a,b \in \Sigma$:
	\be
	\begin{split}
	\sum_{k \geq 0} \alpha_{n,k}(a,a;x) \,z^k
	&=
	\frac{1-z-(n-2)xz}{(1-z+xz)(1-z-(n-1)xz)} 
	\cr
	&=\frac{1}{n}\, \frac{1}{ 1-z-(n-1)\,xz} + \frac{n-1}{n}\, \frac{1}{ 1-z+xz },
	\cr
	\sum_{k \geq 1} \alpha_{n,k}(a,b;x)\, z^k 
	&=
	\frac{zt}{(1-z+xz)(1-z-(n-1)xz)} \cr
	&=\frac{1}{n}\, \frac{1}{1-z-(n-1)\,xz }-\frac{1}{n}\, \frac{1}{1-z+xz},
	\end{split}
	\ee
or equivalently,
	\be\label{eq:alphaformula}
	\begin{split}
	\alpha_{n,k}(a,a;x) &= \frac{1}{n} \left( (1-(n-1)x)^k + (n-1)(1-x)^k\right), \cr
	\alpha_{n,k}(a,b;x) &= \frac{1}{n} \left( (1-(n-1)x)^k -(1-x)^k \right).
	\end{split} 
	\ee
We thus arrive at expressions for the two-point correlation functions:

\begin{prp} \label{lem:sp2pt}
For $a,b \in \Sigma$ with $a \neq b$ and $1 \leq i<j \leq L$,
	\be  \label{twoptform}
	\begin{split}
	\langle \eta^{a,i} \eta^{a,j} \rangle &= \frac 1{n^{2}} + \frac{n-1}{n^{2}}
	\left( \frac{1-x}{1+(n-1)x} \right)^{j-i},\\
	\langle \eta^{a,i} \eta^{b,j} \rangle &= \frac 1{n^{2}} - \frac1{n^{2}}
	\left( \frac{1-x}{1+(n-1)x} \right)^{j-i}.
	\end{split}
	\ee
\end{prp}

\begin{proof}
By Lemma \ref{lem:spindofL} we may assume $i=1$ and $j=L$.
Comparing Lemma \ref{lem:spbarmu} with the definition of the $\alpha_{n,k}(a,b;x)$
in (\ref{eq:defalpha})
we see that for $a,b \in \Sigma$:
	\be
	\langle \eta^{a,1} \eta^{b,L} \rangle = \frac{\alpha_{n,L-1}(a,b;x)}{n(1+(n-1)x)^{L-1}} ,
	\ee  
so that the assertion follows from \ref{eq:alphaformula}.
\end{proof}

The formula \eqref{twoptform} is quite interesting because the first term, 
$1/n^{2}$, has a significance.  From the formula for the density in 
Corollary~\ref{cor:spdens}, we get 
	\be \label{trunc2}
	\langle \eta^{a,1} \eta^{a,L} \rangle - \langle \eta^{a,1}  \rangle \langle  \eta^{a,L} 
	\rangle= \frac{n-1}{n^{2}} \left( \frac{1-x}{1+(n-1)x} \right)^{L-1}.
	\ee
The object on the left hand side is called the {\em truncated} two point correlation function 
in the physics literature, and its value is an indication of  how far the stationary distribution 
is from a product measure. In the case of a product measure, the right hand side would be zero.  
Setting 
	\be
	\alpha=\frac{1-x}{1+(n-1)x},
	\ee
we see that $|\alpha| \leq 1$, and so the truncated correlation function goes exponentially 
to zero as $L \to \infty$. Thus, the stationary measure $\bar\mu_{x}$ behaves like a product 
measure if we do not look for observables which are close to each other.
We can use \eqref{trunc2} to understand one of the differences between the values $x<1$ and $x>1$, 
namely in the way this quantity converges. In the former case, the convergence is monotonic, 
and in the latter, oscillatory.

\section*{Acknowledgements}
The first author (A.A.) would like to acknowledge hospitality 
and support from the Tata Institute of Fundamental Research, Mumbai, India 
where part of this work was done, and thank T. Amdeberhan for discussions.

\bibliographystyle{alpha}
\bibliography{bruijn}

\end{document}